\documentclass[letterpaper, 10 pt, conference]{IEEEconf}  

\IEEEoverridecommandlockouts                              
\usepackage{amsmath,amssymb,graphicx,color,paralist,bm,fancybox}
\newcommand{\eq}{\begin{equation}\begin{array}{lclllllllllllllll}}
\newcommand{\ee}{\end{array}\end{equation}}
\newcommand{\bmt}{\left[ \begin{array}{ccccccccc}}
\newcommand{\emt}{\end{array}\right]}
\newcommand{\bea}{ \begin{eqnarray}}
\newcommand{\eea}{\end{eqnarray}}
\newcommand{\bean}{\begin{eqnarray*}}
\newcommand{\eean}{\end{eqnarray*}}
\newcommand{\bc}{\begin{center}}
\newcommand{\ec}{\end{center}}
\newcommand{\bi}{\begin{itemize}}
\newcommand{\ei}{\end{itemize}}

\title{\LARGE \bf
Al'brekht's Method in Infinite Dimensions
}

\author{Arthur J. Krener
\thanks{This work was supported by AFOSR}
\thanks{A. Krener is with the Department of Applied Mathematics, Naval Postgraduate School,
Monterey, CA 93940, USA
        {\tt\small ajkrener@nps.edu}}%
}%

\begin{document}

\maketitle
\thispagestyle{empty}
\pagestyle{empty}

\begin{abstract}
In 1961 E. G. Albrekht presented a method for the optimal stabilization of smooth, nonlinear, finite dimensional,  continuous time control systems.  This method has been extended to similar systems in   discrete time and to some stochastic systems in continuous and discrete time.  In this paper we extend Albrekht's method to the optimal stabilization of some smooth, nonlinear, infinite dimensional, continuous time control systems whose nonlinearities are described by Fredholm integral operators.
\end{abstract}

\vspace{0.2in}


\section{INTRODUCTION}
A fundamental control engineering problem is to find a feedback law that stabilizes a plant to an operating point.  
Suppose the plant can be modeled by a finite dimensional system of nonlinear differential equations 
\bea \label{nlrdy}
\dot{x}&=& f(z,u)
\eea
and the operating point is $ z=0,u=0$ where $f(0,0)=0$. We assume that the state $z$ is $n$ dimensional
and control $u$ is $m$ dimensional.
We also assume that $ f(z,u)$ is smooth around the operating  point, 
\bean 
f(z,u)&=& Fz+Gu +O(z,u)^2
\eean

Then posing and solving a Linear Quadratic Regulator (LQR) problem will yield a locally stabilizing linear feedback.
We chose an $ (n+m)\times(n+m)$  nonnegative definite matrix
$ [Q,S';S,R]\ge0
$
with $ R>0$  positive definite.
We seek to minimize

\bea \label{lqr}
{1\over 2}\int_0^\infty z'Qz+2z'Su + u'Ru\  dt
\eea
subject to the linear dynamics
\bean \label{lrs}
\dot{z}&=& Fz +Gu
\eean
and a given initial condition $ z(0)=z^0.$

Under the standard assumptions of stabilizability and detectability,  the optimal cost exists and is of the form 
$
{1\over 2} (z^0)'Pz^0
$
and the optimal feedback exists and is of the form 
$
u(t)= Kz(t)
$
where  the $n\times n$ nonnegative definite matrix  $ P\ge 0$  and the $m\times n$ matrix $ K$ satisfy the familiar LQR equations
\eq \label{LQR}
0&=& F'P +PF +Q -(PG+S)R^{-1}(PG+S)'\\
K&=&-R^{-1}(PG+S)' 
\ee
The first equation is called the Algebraic Riccati Equation (ARE).

The function 
$
z'Pz 
$ 
is a local Lyapunov function for the closed loop nonlinear dynamics
using the optimal linear feedback,
\bean
\dot{z}&=& f(z, Kz)
\eean
\bean
{d\over dt} z'(t) Pz(t) &=&((F+GK)z+O(z)^2)'Pz\\
&&+z'P((F+GK)z+O(z)^2)\\
&=& -z'\left( Q+(PG+S)R^{-1}(PG+S)'\right)z\\&&+O(z)^3
\eean

Consider the problem of minimizing a more general criterion
\bea  \label{la}
\int_0^\infty l(z,u) \ dt
\eea
subject to the nonlinear dynamics (\ref{nlrdy})  where the Lagrangian 
is smooth
\bean
l(z,u) &=&{1\over 2}\left( z'Qz+2z'Su + u'Ru\right)+O(z,u)^3
\eean
The higher degree terms in $l(z,u)$ could be penalty terms to ensure
that state and control constraints are satisfied.  They might  destroy  its even symmetry.     The higher degree terms
in $f(z,u)$ might destroy its odd symmetry.

Given that  $z(0)=z^0$ the optimal cost $\pi(z^0) $  if it exists and is smooth  and 
if  the optimal feedback $u(t)=\kappa(z(t))$ exists 
then they satisfy the familiar Hamilton-Jacobi-Bellman (HJB) equations
\eq\label{HJB}
0&=& \mbox{min}_u  \left\{\frac{\partial \pi} {\partial z}(z)f(z,u)+l(z,u)\right\}\\
\\
\kappa(z)&=& \mbox{argmin}_u  \left\{\frac{\partial \pi} {\partial z}(z)f(z,u)+l(z,u)\right\}
\ee
If the quantity to be minimzed in these equations  is a smooth function then the HJB equations 
can be simplified to the sHJB equations
\eq\label{sHJB}
0&=&\frac{\partial \pi} {\partial z}(z)f(z,\kappa(z))+l(z,\kappa(z))\\
\\
0&=& \frac{\partial \pi} {\partial z}(z)\frac{\partial f}{\partial u}(z,\kappa(z))+\frac{\partial f}{\partial u}(z,\kappa(z))
\ee

Assuming $f(z,u)$ and $l(z,u)$ are sufficiently smooth Al'brekht \cite{Al61} showed how to compute the Taylor polynomials of $\pi(z^0) $ and $\kappa(z)$ degree by degree.  At the lowest degrees he obtained the familiar LQR equations.  At higher degrees he obtained a sequence of linear equations for the higher degree coefficients of $\pi(z^0) $ and $\kappa(z)$.  The purpose of this paper is to show that 
Al'brekht's method can be extended to some infinite dimensional control problems.
Navasca   extended Al'brekht's method to discrete time problems \cite{Na02}.  We have extended it to  some stochastic problems in both continuous \cite{Kr19a} and discrete time \cite{Kr19b}.

In the next section we review Al'brekht method for smooth finite dimensional problems.   Rather than rely on the sHJB equations we will 
use a technique which is a conceptually simpler called completing the square.   Completing the square is frequently used  to derive the LQR equations.   We shall use it to find the higher degree terms of the optimal cost.   In Section 4 we extend Al'brekht method to controlled reaction-diffusion systems.   Section 5  contains an example of such a system.

We are not the first to use Al'brekht's method on infinite dimensional systems, see the works of   Kunisch and coauthors \cite{BKP18}, \cite{BKP19}, \cite{KP19}.  Krstic and coauthors have had great success stabilizing infinite dimensional systems where the the nonlinearites are expressed  by Volterra  integral operators of increasing degrees using backstepping techniques, \cite{KS08}, \cite{VK08}.   In our  extension of Al'brekht we assume that the nonlinearites are expressed by Fredholm integral operators of increasing degrees.

\section{Al'brekht Method in Finite Dimension}
Al'brekht assumed that  $ f(z,u)$ and  $ l(z,u)$ are sufficiently smooth to have the Taylor polynomial expansions
\bean
f(z,u)&=& Fz+Gu +f^{[2]}(z,u)+\ldots +f^{[d]}(z,u)\\
&&+O(z,u)^{d+1}\\
l(z,u)&=& {1\over 2}\left(z'Qz+2z'Su + u'Ru\right)+l^{[3]}(z,u)+\ldots \\&&+l^{[d+1]}(z,u)+O(z,u)^{d+2}
\eean
for some degree $ d>1$ where $ ^{[k]}$ indicates terms  of homogeneous degree $ k $ in $ z, u.$

He also assumed that the optimal cost $\pi(z)$ and the optimal feedback have similar 
Taylor polynomial expansions
\bean
\pi(z)&=& {1\over 2}z'Pz+\pi^{[3]}(z)+\ldots +\pi^{[d+1]}(z)+O(z)^{d+2}
\\
\kappa(z)&=& Kz+\kappa^{[2]}(z)+\ldots +\kappa^{[d]}(z)+O(z)^{d+1}\\
\eean

 First we complete the square for the LQR problem.  
 Let $P$ be any symmetric $ n\times n $ matrix.  For any control trjectory $u(t)$ that results in a state trajectory
 $z(t)$ that goes to zero we have
 \small
 \bean
 0&=&{1\over 2}z'(0)Pz(0)+{1\over 2}\int_0^\infty {d\over dt} z'(t)Pz(t)  \ dt \\
 &=& {1\over 2}z'(0)Pz(0)\\
 &&+{1\over 2}\int_0^\infty \bmt z  \\u  \emt' \bmt F'P+PF& PG\\G'P &0\emt \bmt z  \\u  \emt\ dt
 \eean
 \normalsize
 We add this zero quantity to the criterion to be minimized to get a new criterion to be minimized
 \small
 \bean
&& {1\over 2}z'(0)Pz(0)\\
 &&+{1\over 2}\int_0^\infty \bmt z  \\u  \emt' \bmt F'P+PF+Q& PG+S\\G'P+S' &R\emt \bmt z  \\u  \emt\ dt
 \eean
 \normalsize
 We want to choose $P$ and an $ m\times n $ matrix $K$ so the integrand is a perfect square
 \bean
&& \bmt z  \\u  \emt' \bmt F'P+PF+Q& PG+S\\G'P+S' &R\emt \bmt z  \\u  \emt\\
 &&= (u-Kz)'R(u-Kz)
 \eean
 This will be true iff $P$ and $K$ satisfy the LQR equations (\ref{LQR}).  Clearly then $P$ is the kernel of the the optimal cost 
 and $K$ is the optimal feedback gain.
 
 Suppose $\pi^{[3}(z)$ is any homogeneous polynomial of degree three in $z$.
 Again for any control trjectory $u(t)$ that results in a state trajectory
 $z(t)$ that goes to zero we have
 \small
 \bean
 0&=&{1\over 2}z'(0)Pz(0)+{1\over 2}\int_0^\infty {d\over dt} z'(t)Pz(t)  \ dt \\
 &&+ \pi^{[3}(z(0))+\int_0^\infty {d\over dt} \pi^{[3}(z(t)) \ dt
 \\
 &=& {1\over 2}z'(0)Pz(0)\\
 &&+{1\over 2}\int_0^\infty \bmt z  \\u  \emt' \bmt F'P+PF& PG\\G'P &0\emt \bmt z  \\u  \emt\ dt\\
 && + \pi^{[3}(z(0))+\int_0^\infty\frac{\partial\pi^{[3}}{\partial z}(z)\left(F+GK\right)z\\
 && +z'Pf^{[2]}(z,Kz)+l^{[3]}(z,Kz) \ dt +O(z,u)^4
 \eean

 \normalsize
 We add this to the criterion (\ref{la}) to be minimized.  We have already matched quadratic terms.  The cubic terms are
 \bean
&&  \pi^{[3}(z(0))+\int_0^\infty\frac{\partial\pi^{[3}}{\partial z}(z)\left(F+GK\right)z\\
 && +z'Pf^{[2]}(z,Kz) +l^{[3]}(z,Kz)\ dt 
 \eean
 We choose $\pi^{[3}(z)$  so that the integrand vanishes,
\bea
\label{p3}
0=\frac{\partial\pi^{[3]}}{\partial z}(z)\left(F+GK\right)z +z'Pf^{[2]}(z,Kz) +l^{[3]}(z,Kz)
\eea
then $\pi^{[3]}(z^0)$ is the cubic part of the optimal cost.   Notice that the quadratic part
$\kappa^{[2}(z)$ 
of the optimal feedback does not enter in this equation (\ref{p3}).

The solvability of (\ref{p3}) depends on the invertability of the operator
\bea
 \label{op3}
\pi^{[3]}(z)\mapsto \frac{\partial\pi^{[3]}}{\partial z}(z)\left(F+GK\right)z
\eea
acting on   homogeneous polynomials of degree three.  The eigenvalues of this operator
are the sums $\mu_i+\mu_j+\mu_k$ of three eigenvalues of the linear
closed loop dynamics $F+GK$.  Under the standard LQR assumptions these eigenvalues are all
in the open left half plane and hence no triple sum is zero. So the operator (\ref{op3}) is invertible.

Suppose $\psi_1, \ldots, \psi_n $ are the left row eigenvectors of $F+GK$.   Then the corresponding eigenvectors 
of the operator (\ref{op3}) are of the form $(\psi_iz)(\psi_jz)(\psi_kz)$. Let $\Psi$ be the $n\times n $ matrix whose rows are
the $\psi_i$ then the linear change of state coordinates 
\bea 
\zeta&=&\Psi z \label {coc}
\eea
 diagonalizes the operator (\ref{op3})
and makes (\ref{p3}) easy to solve.   

The quadratic terms in  the second sHJB equation (\ref{sHJB}) are
\bean
0&=& \frac{\partial \pi^{[3]}}{\partial z}(z)G+z'P\frac{\partial f^{[2]}}{\partial u}(z,Kz)\\
&&+\frac{\partial l^{[3]}}{\partial u}(z,Kz)+(\kappa^{[2]}(z))'R
\eean
where $\kappa^{[2]}(z)$ is the quadratic part of the optimal feedback.  Since we have assumed $R$ is positive definite
this equation uniquely determines 
$\kappa^{[2]}(z)$.

The higher degrees terms in the optimal cost and optimal feedback are found in a similar fashion.  First we solve 
\bea \label{opk1}
\frac{\partial \pi^{[k+1]}}{\partial z}(z)\left(F+GK\right)z&=&\mbox{Known Stuff}
\eea
and then we solve 
\bean
\kappa^{[k]}(z)=-R^{-1}\left(\frac{\partial \pi^{[k+1]}}{\partial z}(z) G+\mbox{Known Stuff}\right)'
\eean
where the $ \mbox{Known Stuff} $ consist of terms from the Taylor polynomials of $ f(z,u)$ 
and $ l(z,u)$ and previously computed terms from the Taylor polynomials of $ \pi(z)$ 
and $ \kappa(z)$.   The first equation (\ref{opk1}) uniquely determines $  \pi^{[k+1]}(z)$ 
because the eigenvalues of the linear operator 
\bean
\pi^{[k+1]}(z) \mapsto \frac{\partial \pi^{[k+1]}}{\partial z}(z)\left(F+GK\right)z
\eean
are the sums of $k+1$ eigenvalues of $F+GK$.
Again the linear change of state coordinates (\ref{coc})
diagonalizes the equation (\ref{opk1})  for every $ k$.

\section{Controlled Reaction-Diffusion Equations}
In this section we extend Al'brekht's Method to controlled 
reaction-diffusion equations of the form
\bea
\frac{\partial z}{\partial t}(x,t)&=&\frac{\partial^2 z}{\partial x^2}(x,t)+f(x,z(\cdot,t),u(\cdot,t))    \label{crde}\\
\frac{\partial z}{\partial x}(0,t)=&0&=\frac{\partial z}{\partial x}(1,t)\label{noflux}\\
z(x,0)&=& z^0(x)
\eea
for $ x\in [0,1]$  with Neumann (no flux) boundary conditions.    We wish to stabilize $ z(x,t)$ to $ z(x)=0$. 
Both $ z(x,t)$  and $ u(x,t)$  could be vector valued but for simplicity of exposition we assume that they are scalar valued.

Notice that the reaction term is a functional of $ z(\cdot,t),u(\cdot,t)$. 
We assume that  $ f(x,z(\cdot,t),u(\cdot,t)) $ is given by 
a sum of Fredholm integral operators of increasing degrees
\bean \label{f}
&& f(x,z(\cdot),u(\cdot)) \\
&&=\int_0^1 F^{[1]}(x,x_1)z(x_1) +G^{[1]}(x,x_1)u(x_1)\ dx_1 \nonumber \\
&&+\int_0^1\int_0^1 F^{[2]}(x,x_1,x_2)z(x_1)z(x_2)\ dx_1dx_2
\\ &&+\ldots \nonumber
\eean 

This is not the most general  reaction term that we could consider, for example we could have terms quadratic in $ u(\cdot)$  or bilinear in $ z(\cdot)$ and $ u(\cdot)$ or higher degrees
in $ z(\cdot)$ and $ u(\cdot)$. 
WLOG we assume all the Fredholm kernels are symmetric with respect to the subscripted $ x_i,$ e.g.
\bean
F^{[2]}(x,x_1,x_2)&=&F^{[2]}(x,x_2,x_1)
\eean
Note that $ F^{[k]}(x,x_1,\ldots, x_k)$   and $ G^{[1]}(x,x_1)$  could be generalized functions.
For example if $ F^{[1]}(x,x_1)=F(x)\delta(x-x_1)$  and $ G^{[1]}(x,x_1)=G(x)\delta(x-x_1)$  then the linear part of the reaction term  is
\bea \nonumber
&&\int_0^1 F^{[1]}(x,x_1)z(x_1) +G^{[1]}(x,x_1)u(x_1)\ dx_1\\
&&= F(x)z(x) +G(x)u(x) \label{linear}
\eea
To keep the notation  simple we shall assume that (\ref{linear}) is the case. 

To find a linear  feedback that locally stabilizes the  system to $ z(x)=0$,
we pose the infinite dimensional LQR problem of
minimizing
\bea \label{crt}
&&{1\over 2}\int_0^\infty \int_0^1 |z(x,t)|^2+|u(x,t)|^2 \ dx\  dt\\&=& \nonumber
{1\over 2}\int_0^\infty \int_0^1 \int_0^1 \delta(x_1-x_2)z(x_1,t)z(x_2,t)\\
&&+\delta(x_1-x_2)u(x_1,t)u(x_2,t) \ dx_1dx_2\  dt \nonumber
\eea
subject to the linear part of the dynamics.  We pose this simple Lagrangian  but our method readily extends
to more complicated ones in Fredholm form..

We complete the square again.  Suppose we have a Fredholm quadratic  form  in $z^0(x)$
\bean
\int_0^1\int_0^1 P^{[2]}(x_1,x_2)z^0(x_1)z^0(x_2)\ dx_1 dx_2
\eean
and the control trajectory $u(x,t)$  takes $z(x,t)$ to zero as $t\to \infty$. Then
\bean 
0= \int_0^1\int_0^1 P^{[2]}(x_1,x_2)z^0(x_1)z^0(x_2)\ dx_1 dx_2\hspace{.7in}\\
+ \int_0^\infty \int_0^1\int_0^1 {d\over dt} \left(P^{[2]}(x_1,x_2)z(x_1,t)z(x_2,t)\right)\ dx_1 dx_2 dt
\eean

We add this zero quantity  to the criterion (\ref{crt}) to get a new expression
to be minimized 
\eq \nonumber
\int_0^1 \int_0^1 P^{[2]}(x_1,x_2) z^0(x_1)z^0(x_2) \ dx_1 dx_2\\
\\
 + \int_0^\infty\int_0^1 \int_0^1 P^{[2]}(x_1,x_2)\\
 \\
\times \left(\frac{\partial^2 z}{\partial x_1^2}(x_1,t)z(x_2,t)+ z(x_1,t)\frac{\partial^2 z}{\partial x_2^2}(x_2,t)
\right) \\
\\
 +P^{[2]}(x_1,x_2)\left(F(x_1)z(x_1,t)+G(x_1) u(x_1,t)\right)z(x_2,t)\\
 \\
 +P^{[2]}(x_1,x_2)z(x_1,t)\left(F(x_2)z(x_2,t)+G(x_2)u(x_2,t)\right)\\
 \\
+ \delta(x_1-x_2)z(x_1,t)z(x_2,t)\\
\\
+\delta(x_1-x_2)u(x_1,t)u(x_2,t) \ dx_1 dx_2dt 
\ee

We have assumed that $z(x,t)$ satisfies Neumann boundary conditions.  If we  assume that $P^{[2]}(x_1,x_2)$  also satisfies Neumann boundary conditions then  when we integrate by parts twice we get 
\eq \nonumber
\int_0^1 \int_0^1 P^{[2]}(x_1,x_2) z^0(x_1)z^0(x_2) \ dx_1 dx_2\\
\\
+ \int_0^\infty\int_0^1 \int_0^1 \left(\Delta P^{[2]}(x_1,x_2)+\delta(x_1-x_2)\right)z(x_1,t)z(x_2,t)\\
\\
 +P^{[2]}(x_1,x_2)\left(F(x_1)z(x_1,t)+G(x_1) u(x_1,t)\right)z(x_2,t)\\
 \\
 +P^{[2]}(x_1,x_2)z(x_1,t)\left(F(x_2)z(x_2,t)+G(x_2)u(x_2,t)\right)\\
 \\
+\delta(x_1-x_2)u(x_1,t)u(x_2,t) \   dt  \ dx_1 dx_2dt 
\ee

\normalsize
where $
\Delta P^{[2]}(x_1,x_2)$ is the two dimensional Laplacian of  $ P^{[2]}(x_1,x_2)$.

We want the integrand of the time integral to be a perfect square of the form
\eq  \nonumber
\int_0^1 \int_0^1 \left(u(x_1,t)-\int_0^1K^{[1]}(x_1,x_3)z(x_3,t)\ dx_3\right)  \delta(x_1-x_2)\\
\times \left(u(x_2,t)-\int_0^1K^{[1]}(x_2,x_3)z(x_3,t)\ dx_3\right)\ dx_1dx_2
\ee
for some not necessarily symmetric $K^{[1]}(x,x_3)$.    

This leads to infinite dimensional LQR equations 
\bea
&&K^{[1]}(x,x_1)= -P^{[2]}(x,x_1)G(x_1)\label{K1}\\
&&\int_0^1 P^{[2]}(x_1,x_3)G(x_3)G(x_3)P^{[2]}(x_3,x_2)\ dx_3\nonumber\\
&&= \Delta P^{[2]}(x_1,x_2)+\delta(x_1-x_2) \nonumber\\
&&+
F(x_1)P^{[2]}(x_1,x_2) +P^{[2]}(x_1,x_2)F(x_2) \label{P2}
\eea
The second of these equation is called a Riccati PDE and it is to be interperted 
in the weak sense.  If $\theta(x_1),\theta(x_2) $ are any $C^2$ functions then
\eq \nonumber
\int_0^1\int_0^1\int_0^1\theta(x_1) P^{[2]}(x_1,x_3)G(x_3)\\
\\
 \times G(x_3)P^{[2]}(x_3,x_2)\ \theta(x_2)\ dx_1dx_2dx_3\\
\\
=\int_0^1\int_0^1\theta(x_1) \left(\Delta P^{[2]}(x_1,x_2)+\delta(x_1-x_2)\right. \\
\\
\left.+
F(x_1)P^{[2]}(x_1,x_2)+P^{[2]}(x_1,x_2)F(x_2)\right)\theta(x_2)\ dx_1dx_2
\ee
Similar equations have appeared in the works
of J.L. Lions \cite{JL71},  J. Burns \cite{BRK94}, \cite{BK95},  K. Hulsing \cite{Hu99}, \cite{BH01}, B. Batten King \cite{Ki00} and others.

If $P^{[2]}(x_1,x_2)$ is a weak solution of the Riccati PDE (\ref{P2}) then it is the kernel of the degree two Fredholm form that is the quadratic part of the optimal cost.   The optimal linear feedback gain is given by (\ref{K1}).
The closed loop linear dynamics is
\bean
\frac{\partial z}{\partial t}(x,t)&=&\frac{\partial^2 z}{\partial x^2}(x,t)+F(x)z(x,t)\\
&&+\int_0^1 G(x)K^{[1]}(x,x_1)z(x_1,t)\ dx_1
\eean

A standard approach to solving the Riccati PDE (\ref{P2})  is to expand $ P^{[2]}(x_1,x_2)$ in the eigenfunctions of the diffusion operator.   The eigenvectors and the eigenvalues of ${\partial^2\over \partial x^2}$ subject to Neumann boundary conditions on $x\in[0,1]$ are
\eq \label{oevv}
\lambda_0=0,&&\phi_0(x)=1\\
\lambda_i=-i^2\pi^2,&&\phi_i(x)=\sqrt{2} \cos(i\pi x)
\ee
for $i=0,1,2,\ldots$.  Notice that this is an orthonormal family
If we assume that 
\bea \label{PP2}
P^{[2]}(x_1,x_2)&=& \sum_{i,j=0}^\infty \Pi_{i,j}\phi_i(x_1)\phi_j(x_2)
\eea
Then the Riccati PDE (\ref{P2}) becomes an Algebraic Riccati equation for the
 infinite dimensional matrix $[\Pi_{i,j}]$.

Let the linear closed loop eigenvalues and left eigenvectors be denoted by
\bean
\mu_i,&& \psi_i(x)
\eean
With the criterion (\ref{crt}) the LQR is clearly detectable.  Under  an additional assumption of stabilizability we have that
all the $\mu_i$ are in the open left half plane \cite{CZ95}. But in general the $\psi_i(x)$ are not orthonormal.

Let $P^{[3]}(x_1,x_2,x_3)$ be the kernel of any degree three Fredholm form and suppose
that $ u(x,t) $ is a control trajectory that results in an asymptotically  stable  state trajectory
$ z(x,t). $  Then
\small
\eq \nonumber
0={1\over 2}\int_0^1\int_0^1 P^{[2]}(x_1,x_2)z^0(x_1)z^0(x_2) \ dx_1 dx_2\\
\\
+{1\over 2}
\int_0^\infty \int_0^1\int_0^1 {d\over dt}\left(P^{[2]}(x_1,x_2)z(x_1,t)z(x_2,t) \right)\ dx_1 dx_2dt\\
\\+ \int_0^1\int_0^1\int_0^1 P^{[3]}(x_1,x_2,x_3)z^0(x_1)z^0(x_2)z^0(x_3) \ dx_1 dx_2 dx_3\\
\\
+ 
\int_0^\infty \int_0^1\int_0^1\int_0^1 {d\over dt}\left(P^{[3]}(x_1,x_2,x_3)z(x_1,t)z(x_2,t)z(x_3,t) \right)\\
\\
 dx_1 dx_2dx_3 dt
\ee

\normalsize
As before we add this zero quantity to the criterion to be minimized to get a new expression to be minimized.  
We have already matched up the quadratic terms in this expression.
The cubic terms   are 
\small
\eq
\nonumber
\int_0^1\int_0^1\int_0^1 P^{[3]}(x_1,x_2,x_3)z^0(x_1)z^0(x_2)z^0(x_3) \ dx_1 dx_2 dx_3\\
\\
 + 
\int_0^\infty \int_0^1\int_0^1\int_0^1 \left( \Delta P^{[3]}(x_1,x_2,x_3) \right.\\
\\
+3\int_0^1 P^{[3]}(x_4,x_1,x_2)\left(F(x_4)+G(x_4)K^{[1]}(x_4,x_3)\ dx_4 \right) 
\\
\\ 
\left.
+\int_0^1P^{[2]}(x_4,x_1)F^{[2]}(x_4,x_2,x_3)\ dx_4\right)\\
\\
\times \ z(x_1,t) z(x_2,t)z(x_3,t)\ dx_1 dx_2 dx_3 dt
\ee
where $\Delta P^{[3]}(x_1,x_2,x_3)$ is the three dimensional Laplacian.

We set the  time integrand  to zero to get a weak linear elliptic PDE for the symmetric function $ P^{[3]}(x_1,x_2,x_3)$
\eq
\label{P3}
0=\Delta P^{[3]}(x_1,x_2,x_3)\\
\\+3\int_0^1 P^{[3]}(x_4,x_1,x_2)\left(F(x_4)+G(x_4)K^{[1]}(x_4,x_3)\ dx_4 \right) \\ 
\\+ \int_0^1P^{[2]}(x_4,x_1)
F^{[2]}(x_4,x_2,x_3)\ dx_4
\ee
subject to Neumann boundary conditions.   By weak we mean if we multiply the right side of  this equation by any three $C^2$ functions
 $\theta(x_1)\theta(x_2)\theta(x_3)$ and integrate  over the unit cube we get zero.

Based on what we saw for finite dimensional systems it is natural to make the linear change of state coordinates
\bean
\zeta_i(t)&=& \int_0^1 \psi_i(x)z(x,t)
\eean
Suppose
\bean
P^{[3]}(x_1,x_2,x_3)&=& \sum_{i,j,k=0}^\infty\Pi_{i,j,k}\psi_i(x_1)\psi_j(x_2)\psi_k(x_3)
\eean
Then (\ref{P3}) becomes the triple sequence of equations
\eq \label{P3a}
0=(\mu_i+\mu_j+\mu_k)\ \Pi_{i,j,k} \\
\\
+\int_0^1\int_0^1\int_0^1  \int_0^1P^{[2]}(x_4,x_1)
F^{[2]}(x_4,x_2,x_3)\\
\\
\times \psi_i(x_1)\psi_j(x_2)\psi_k(x_3)\ dx_1dx_2dx_3 dx_4
\ee
This determines $\Pi_{i,j,k} $.

The quadratic part of the optimal feedback is then found from the second 
sHJB equation,
\bean
K^{[2]}(x,x_1,x_2)&=& -3P^{[3]}(x,x_1,x_2)G
\eean
The higher degree terms are found in a similar fashion

 \begin{figure}
  \label{fpcf}
 \centering
\includegraphics[width=3in]{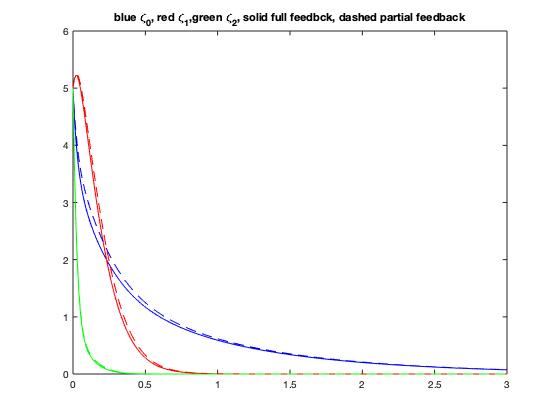}
\caption{First three modes under full (solid) and partial (dashed) cubic feedback}     
 \end{figure}

\section{Example}

We close with a simple example that is a quadratic modification of Example 6.2   of Curtain and Zwart \cite{CZ95}.
Consider a rod of length one with distributed heating/cooling, no flux 
boundary conditions and a quadratic nonlinearity.
\bean
\frac{\partial z}{\partial t}(x,t)&=&\frac{\partial^2 z}{\partial x^2}(x,t)+u(x,t)+z(x,t)^2    \label{crde}\\
\frac{\partial z}{\partial x}(0,t)=&0&=\frac{\partial z}{\partial x}(1,t)\label{noflux}\\
z(x,0)&=& z^0(x)
\eean
Then
\bean
F(x)=0,&&
G(x)= 1\\
F^{[2]}(x,x_1,x_2)&=& \delta(x-x_1)\delta(x-x_2)
\eean
and all the other higher degree terms are zero.

To find a feedback that stabilizes this system to $z(x)=0$
we pose the optimal control problem of minimizing (\ref{crt}) subject 
to this dynamics.  

The open loop linear eigenvalues
and eigenvectors are (\ref{oevv}) so if we assume that (\ref{PP2}) then
the Riccati PDE becomes 
\bean
\sum_{k=0}^\infty \Pi_{i,k}\Pi_{k,j}&=& (\lambda_i+\lambda_j) \Pi_{i,j}+\delta_{i,j}
\eean
This is identical to equation (6.81) of \cite{CZ95} and their solution is
\bean
\Pi_{i,j}&=& \delta_{i,j}\left(\lambda_i+\sqrt{\lambda^2_i+1}\right)\\
\eean 
The linear part of the optimal feedback is 
\bean 
u(x,t)&=& \int_0^1 K^{[1]}(x,x_1)z(x_1,t)\ dx_1\\&=&-\int_0^1 P^{[2]}(x,x_1)z(x_1,t)\ dx_1
\eean
The closed loop eigenvalues are
$
\mu_i=-\sqrt{\lambda_i^2+1}
$.
The  closed loop
eigenvectors  are still
$
\phi_i(x)
$.
But the basin of asymptotic stability of the nonlinear system closed by the 
linear feedback is not very large. If $ z^0(x)=1+\epsilon  $ for $ \epsilon >0$ then the linearly closed loop 
nonlinear system diverges.

The Fredholm kernel of cubic part of the optimal cost
is of the form 
\bean
P^{[3]}(x_1,x_2,x_3)&=&\sum_{i,j,k=0}^\infty \Pi_{i,j,k}\ \phi_i(x_1)\phi_j(x_2)\phi_k(x_3)
\eean
and (\ref{P3a}) becomes 
\bean
0&=& (\mu_i+\mu_j+\mu_k)\Pi_{i,j,k}+{1\over 2}\left(\Pi_{i,j+k}+\Pi_{i,j-k}\right)
\eean 
The solutions $\Pi_{i,j,k}$ of these equations are not symmetric in the indices $i,j,k$.  They 
need to be symmetrized.
The Fredholm kernel of the quadratic part of the optimal feedback is given by
\bean
K^{[2]}(x,x_1,x_2)&=&-3P^{[3]}(x,x_1,x_2)
\eean

The Fredholm kernel of quartic part of the optimal cost
is of the form 
\bean
P^{[4]}(x_1,x_2,x_3,x_4)= \sum_{i,j,k,l=0}^\infty \Pi_{i,j,k,l}\ \phi_i(x_1)\phi_j(x_2)\phi_k(x_3) \phi_l(x_4)
\eean
where the coefficients satisfy the equations
\bean
0&=& \left(\mu_i+\mu_j+\mu_k+\mu_l\right)\Pi_{i,j,k,l}+{3\over 2}\left( \Pi_{i,j,k+l}+\Pi_{i,j,k-l} \right)\\
&& -{9\over 2}\sum_{r=0}^\infty \Pi_{i,j,r}\Pi_{r,k,l}
\eean
Again the solutions $\Pi_{i,j,k,l}$ of these equations are not symmetric in the indices $i,j,k,l$.  They 
need to be symmetrized.
The Fredholm kernel of the cubic part of the optimal feedback is given by
\bean
K^{[3]}(x,x_1,x_2,x_3)&=&-4P^{[3]}(x,x_1,x_2,x_3)
\eean

The linear closed loop poles  go to $ -\infty$  quite fast, for example $ \mu_3\approx -88.8321$ so we did a Galerkin projection on the first
three eigenfunctions.
Let $ \zeta_i(t)=<\phi_i(x), z(x,t)>$ and $ \nu_i(t)=<\phi_i(x), u(x,t)>$   for $ i=0,1,2.$  
The Galerkin projection
is the three dimensional nonlinear control system
\bean
\dot{\zeta}_0&=& \lambda_0 \zeta_0+\nu_0+ \zeta_0^2+{1\over2}\zeta_1^2+{1\over 2}\zeta_2^2\\
\dot{\zeta}_1&=&\lambda_1 \zeta_1+\nu_1+2\zeta_0\zeta_1+\zeta_1\zeta_2\\
\dot{\zeta}_2&=&\lambda_2\zeta_2+\nu_2+2\zeta_0\zeta_2+{1\over 2}\zeta_1^2
\eean

Using our Nonlinear Systems Toolbox we found  
the optimal  cost   $ \pi(\zeta)$  to degree $ 4$ for the three dimensional Galerkin approximation.
If the initial state   is $ [\zeta_0; \zeta_1;\zeta_2]$ then

\bean
\pi^{[2]}( \zeta)&=&0.5000 \ \zeta_0^2+   0.0253 \ \zeta_1^2+0.0063\ \zeta_2^2\\
\pi^{[3]}( \zeta)&=&0.3333    \ \zeta_0^3  +   0.0288  \ \zeta_0 \zeta_1^2  
\\&& +0.0066  \ \zeta_0\zeta_2^2
+    0.0010\ \zeta_1^2 \zeta_2  \\
\pi^{[4]}( \zeta)&=& 0.1250 \ \zeta_0^4+   0.0281 \ \zeta_0^2 \zeta_1^2+        0.0065  \ \zeta_0^2 \zeta_2^2\\&&    +       0.0012 \ \zeta_0 \zeta_1^2\ \zeta_2     + 0.0004 \ \zeta_1^4  \\&&    +     0.0002 \  \zeta^2_1\zeta_2^2+0.0000\ \zeta_2^4
\eean

By way of comparison 
if the initial state  of the infinite dimensional system is
$
z^0(x)=\sum_{i=0}^\infty \zeta_i \phi_i(x)
$  
and if $\Pi^{[k]}( \zeta)$ is degree $k$ part of its optimal cost then
\bean
\Pi^{[2]}( \zeta)&=&0.5000 \ \zeta_0^2+0.0253 \ \zeta_1^2 +0.0063\ \zeta_2^2\\
&&+{1\over 2} \sum_{i=3}^\infty \Pi_{i,i} \zeta_i^2
\\
\Pi^{[3]}( \zeta)&=&0.3333\ \zeta_0^3 + 0.0276 \ \zeta_0\zeta_1^2\\
&&+0.0065\ \zeta_0\zeta_2^2+0.0010\ \zeta_1^2\zeta_2\\
&&+\sum_{i,j,k=3}^\infty\Pi_{i,j,k}\zeta_i\zeta_j\zeta_k\\
\Pi^{[4]}( \zeta)&=&0.1250 \ \zeta_0^4+   0.0279 \ \zeta_0^2 \zeta_1^2+        0.0063  \ \zeta_0^2 \zeta_2^2\\&&+           0.0011 \zeta_0 \zeta_1^2\ \zeta_2+ 0.0004 \ \zeta_1^4 \\&& +     0.0002 \  \zeta^2_1\zeta_2^2
+0.0000\ \zeta_2^4
\\
&&+\sum_{i,j,k,l=3}^\infty\Pi_{i,j,k,l}\zeta_i\zeta_j\zeta_k\zeta_l
\eean

We started the three dimensional Galerkin approximation 
at $ \zeta_0=5,\zeta_1=5,\zeta_2=5 $ which is way outside the basin of asymptotic
stability for the linearly closed loop nonlinear system
 We used two cubic  feedbacks.
The  full cubic feedback was the one computed by the Nonlinear Systems Toolbox \cite{NST} using Al'brekht's method on the three dimensional system.   But under the linear closed loop dynamics the  monomial $ \zeta_{i_1}\cdots \zeta_{i_k}$
satisfies
\bean
{d\over dt}\zeta_{i_1}\cdots \zeta_{i_k}&=&\left(\mu_{i_1}+\cdots +\mu_{i_k}\right) \zeta_{i_1}\cdots \zeta_{i_k}
\eean
so it may go to zero extremely fast.   
The partial cubic feedback ignores any monomial $ \zeta_{i_1}\cdots \zeta_{i_k}$  where 
\bean
\mu_{i_1}+\cdots +\mu_{i_k} &<& -20
\eean
and performs almost as well as the full cubic feedback, Figure 1. 

 Because the linear change of coordinates (\ref{coc}) diagonalizes the 
equations (\ref{opk1}) for every $k\ge 2$ we don't have to solve these equations for all the monomials of degree $k+1$.  We can choose to solve them only for the monomials that don't go to zero very fast under the linear closed loop dynamics

\section{Conclusion}  Al'brekht developed a method to compute the Taylor polynomial expansions 
for the optimal cost and optimal feedback for smooth, finite dimensional,   infinite horizon optimal control problems whose
linear quadratic part satisfies the standard LQR conditions.  We have shown that Al'brekht's method can be extended
to the optimal stabilization  of some smooth infinite dimensional controlled reaction-diffusion equations whose
linear quadratic part satisfies the infinite dimensional  LQR conditions.   The crucial steps in the extension are the ability
to express the nonlinearities as Fredholm linear operators and the ability to compute a few of the least   stable left eigenfunctions of closed loop linear part of the dynamics.

  The example that we presented had distributed control.  The next step is to apply Al'brekht's merhod to a  reaction-diffusion problem 
under boundary control.  We would also like to extend Al'brekht's method to other infinite dimensional problems such as nonlinear wave equations and nonlinear delay equations.

    \end{document}